\documentclass[10pt,a4paper,leqno,openright]{amsart}

\date{}

\usepackage{latexsym}

\usepackage{amsfonts}

\usepackage{amsmath}

\usepackage{amssymb}

\usepackage{amscd}

\input xy
\xyoption{all}

\newtheorem{theorem}{Theorem}[section]
\newtheorem{lemma}{Lemma}[section]

\newtheorem{claim}{Claim}

\newcommand{\ma}[1]{\ensuremath{\mathbb{#1}}}
 
\font\bb=msbm11 at 11 pt

\def \Z {\hbox{\bb Z}}
\def \Q {\hbox{\bb Q}}
\def \N {\hbox{\bb N}}
\def \P {\hbox{\bb P}}
\def \F {\hbox{\bb F}}

\def \A {\hbox{\bb A}}
\def \Sp {\mbox{\rm{Spec}}}
\def \Ho {\mbox{\rm{Hom}}}

\def \T {\mathcal{T}}
\def \FF {\mathcal{F}}

\def \RR {\mathcal{R}}

\def \n {\mathcal{N}}
\def \L {\mathcal{L}}
\def \O {\mathcal{O}}

\def \Q {\hbox{\bb Q}}

\newcommand{\Spec}{\ensuremath{\mbox{\rm{Spec }}}}
\newcommand{\Proj}{\ensuremath{\mbox{\rm{Proj }}}}
\newcommand{\Gal}{\ensuremath{\mbox{\rm{Gal }}}}
\newcommand{\rp}{\ensuremath{\mbox{\rm{rp}}}}
\newcommand{\Hom}{\ensuremath{\mbox{\rm{Hom}}}}
\newcommand{\Ker}{\ensuremath{\mbox{\rm{Ker}}}}
\newcommand{\Res}{\ensuremath{\mbox{\rm{Res}}}}

\newcommand{\Tr}{\ensuremath{\mbox{\rm{Tr}}}}
\newcommand{\Id}{\ensuremath{\mbox{\rm{Id}}}}

\title{Lifts of points on curves and exponential sums. }
\author{R\'{e}gis Blache}

\begin{document}

\maketitle
 
\centerline{\'Equipe ``Alg\`ebre Arithm\'etique et Applications''}
\centerline{Universit\'e Antilles Guyane}
\centerline{Campus de Fouillole}
\centerline{97159 Pointe \`a Pitre CEDEX - FWI}
\centerline{Email :  {\tt rblache@univ-ag.fr}}

\bigskip

\begin{abstract}
We give bounds for exponential sums over curves defined over Galois rings. We first define summation subsets as the images of lifts of points from affine opens of the reduced curve, and give bounds for the degrees of their coordinate functions. Then we get bounds for exponential sums, extending results of
Kumar et al., Winnie Li over the projective line, and Voloch-Walker over elliptic curves and $C_{ab}$ curves.
\end{abstract}
\medskip


\section{Introduction}

Character sums estimates over finite fields have been widely studied since Gauss. During the last decade, coding theorists initiated the study of a new type of sums, over points with coordinates in a $p$-adic field. The first result in this direction was the generalisation of the Weil-Carlitz-Uchiyama bound. Let $f\in \O_m[X]$ be a polynomial over the ring of integers of $K_m$, the unramified extension of degree $m$ of the field of $p$-adic numbers $\Q_p$, $\T=\T^*\cup\{0\}$ the Teichm\"uller of $\O_m$, where $\T^*$ is the multiplicative subgroup of elements of order prime to $p$ in $\O_m^*$ ; if we write the $p$-adic expansion of $f$ as $f_0+pf_1+\dots+p^{l-1}f_{l-1}+\dots$, $f_i\in\T[X]$, Kumar et al. (cf. \cite{khc}) gave a bound for the following exponential sum, where $\psi$ is an additive character of order $p^l$ of $\O_m$ :
$$\left| \sum_{x\in \T} \psi(f(x))\right| \leq \left(\max_{0\leq i\leq l-1}(p^{l-1-i}\deg f_i)-1\right)p^{\frac{m}{2}}.$$
For convenience, we choose to work over ``Galois rings'', i.e. the rings $R:=\O_m/p^l\O_m:=GR(p^l,m)$. Note that $\T$ is the image of a section of reduction modulo $p$ : $R \rightarrow \F_{p^m}$. In other words, it is a ``lift of points'' from $\A_{\F_{p^m}}^1(\F_{p^m})$ to $\A_{R}^1(R)$. The purpose of this paper is to give bounds for analogous sums, defined on any smooth projective curve over $R$, not only the projective line. 

\medskip

To do so, we adopt the following point of view : let $w: R \rightarrow W_l(k)$ be the isomorphism between the Galois ring $R$ and the ring of Witt vectors of length $l$ with coefficients in $k:=\F_{p^m}$. If $C$ is a smooth projective curve over $R$, equiped with an embedding $C\subset \P_R^n$, let $U=\Sp B\subset \A_R^n=\Sp R[X_1,\dots,X_n]$ be an affine open of $C$ of the form $C\backslash \{C\cap H\}$, with $H$ an hyperplane in $\P_R^n$ ; its image by Greenberg's functor is, roughly speaking, the subvariety $\FF U\subset \A_k^{nl}$ whose points are the $(w(x_1),\dots,w(x_n))\in k^{nl}$, where $(x_1,\dots,x_n)$ describe the points of $U(R)$. We can define a ``reduction morphism'' $\rho_U:\FF U \rightarrow U_k:=U\otimes k=\Sp B_0$ sending $(w(x_1),\dots,w(x_n))$ to $(\overline{x_1},\dots,\overline{x_n})$, i.e. each $w(x_i)$ on its first coordinate, which is the image of $x_i$ by reduction modulo $p$ $R\rightarrow k$. Since $U$ is smooth, there exists a section $\sigma$ of $\rho_U$ and a morphism of $R$-algebras $\Gamma(s):B\rightarrow W_l(B_0)$ associated to $\sigma$. Then $\sigma(U_k)(k)\subset \FF U(k)$ corresponds to a subset $\T_\sigma\subset U(R)$, and if $f\in B$ is a regular function on $U$, we can rewrite (theorem 1.2) the character sum :
$$\sum_{x\in \T_\sigma} \psi(f(x)) = \sum_{x\in U_k(k)} \psi(\Gamma(s)(f)(x)).$$
To give a bound for this last sum, we use Artin-Schreier-Witt theory (the theory of abelian extensions of a field of characteristic $p$, whose exponent is a power of $p$) ; it allows us to associate an idele class character to a Witt vector of functions. From Riemann hypothesis for curves over a finite field (Weil's theorem), giving a bound for the last character sum reduces to the calculation of the conductor of such an idele class character, which depends on the pole order of the components of the Witt vector of functions $\Gamma(s)(f)$.

\medskip

Let us see how this works if $C$ is the projective line $\P_R^1$, and $U$ the affine line $\A_R^1$. In this case, $\FF U$ is affine $l$-space $\A_k^l$, and the reduction morphism is projection on the first coordinate $\A_k^l \rightarrow \A_k^1$. Now since $w(\T)$ is the set $\{(x,0,\dots,0),~x\in k\}$, $\T$ corresponds to image of the section $\sigma :\A_k^1 \rightarrow \A_k^l$ defined by $\sigma(x)=(x,0,\dots,0)$. The morphism $\Gamma(s) :R[X]\rightarrow W_l(k[X])$ is defined by $\Gamma(s)_{|R}=w$, and $\Gamma(s)(X)=(X,0,\dots,0)$. Finally, if $f=f_0+\dots+p^{l-1}f_{l-1}$ is as above, a Witt vector calculation shows that the conductor of the idele class character associated to $\Gamma(s)(f)$ is at most $\max_i(p^{l-1-i}\deg f_i)+1$, and we get the bound.    

\medskip

In the general case, the section $\sigma$ is much harder to write down, and many efforts are devoted to make it as explicit as possible. Actually we give a bound for the pole orders of the ``coordinate functions'' of certain sections $\sigma$, i.e. the components $X_{i0}, X_{i1},\dots,X_{il-1}$ of the vectors $\Gamma(s)(X_i)$ : we show (theorem 2.1) that if $C_k=C\otimes k$ has genus $g$, and if $D_0$ is the intersection divisor of $C_k$ with $H_0\otimes k$ in $\P_k^n$, then there exists a section with $X_{ij}$ having polar divisor less than $p^j(jD+(j+1)D_0)$, where $D$ is a divisor of degree $2g-2+\lceil \frac{2g-1}{p}\rceil$ with support contained in that of $D_0$. Similar calculations have been made in the case of elliptic and hyperelliptic curves by Finotti (cf \cite{fin}), and this result may be of independent interest from the rest of the article.

\medskip

Let us say a few words about the bounds we get, and compare them with previously known results. As in Weil-Carlitz-Uchiyama bound, the conductor is a ``weighted degree'' of the pole orders of the functions in the $p$-adic expansion of the function, inflated by a term of order $2gp^{l-1}$ depending only on the genus of the reduced curve and on the ring $R$. This generalises known results : the case of projective line had been treated by Winnie Li (cf. \cite{li}), who gives the exact degree of the conductor of the idele class character associated to additive or multiplicative characters and rational functions. Our method gives the same result in the case of additive characters (cf remark at the end of section 4) ; this shows that the bounds we give are tight with respect to the weighted degree of the pole orders. Other examples of such character sums in the case $l=2$ can be found in \cite{vw1}, \cite{vw2}, respectively for elliptic and $C_{ab}$ curves ; these results are particular cases of the ones we obtain. They are used to construct error correcting codes over Galois rings with good parameters. Note also that the Weil Carlitz Uchiyama bound allows the construction of large families of periodic sequences with low cross correlation.

\medskip

The paper is organized as follows : in section 1 we show the existence of the sections described above, and show how they allow us to rewrite character sums (theorem 1.2) ; for certain of these sections, we give bounds for the pole orders of their coordinate functions in theorem 2.1 in section 2. Section 3 is independent from the rest of the article ; we recall without proof some results of Artin Schreier Witt theory ; the main result is theorem 3.1, which gives bounds for exponential sums associated to a Witt vector of functions and an additive character of order a power of $p$. We suspect this result to be well known, but we din't find any reference (theorem 3.1 in \cite{vw2} covers only the case of functions having a single pole at a rational point). It reduces our bounds to the calculation of the reduced pole orders of certain Witt vectors, and this is done in section 4. In the last section, we give the bounds, first in the general case (theorem 5.1), then we give two corollaries in particular cases where the bounds are simpler.
   
\bigskip

\centerline{{\sc Notations}}

\medskip

We first recall notations about rings of Witt vectors ; let $A$ be a ring of characteristic $p$ : in the sequel, $W_l(A)$ denotes the ring of Witt vectors of length $l$ with coefficients in the ring $A$. As usual, 
$$V:W_l(A)\rightarrow W_l(A),~V(a_0,\dots,a_{l-1})=(0,a_0,\dots,a_{l-2})$$ 
denotes the {\it Verschiebung}, and :
$$F:W_l(A)\rightarrow W_l(A),~F(a_0,\dots,a_{l-1})=(a_0^p,\dots,a_{l-1}^p)$$ 
the Frobenius of $W_l(A)$. We recall the following relations, which we use in the sequel :
$$ FVx=VFx=px~;~ V^axV^by=V^{a+b}(F^bxF^ay)~\text{for all}~x,y\in W_l(A),~a,b \in \N.$$

Let $R=GR(p^l,m)$ be a Galois ring with residue field $k=\F_{p^m}$, and $R_d=GR(p^l,md)$ its unramified extension of degree $d$, with residue field $k_d=\F_{p^{md}}$. Recall that there is a canonical isomorphism $w:R_d\rightarrow W_l(k_d)$. Recall that $R_d^*$ has a subgroup of order $p^{md}-1$ : we denote it by $\T_d^*$ ; we also let $\T_d=\T_d^*\cup \{0\}$ and call it the ``Teichm\"uller'' of $R_d$. Note that via $w$, $\T_d$ corresponds to the elements of the form $(a,0,\dots,0)$, $a\in k_d$.

\medskip

We will need the following easy lemma about valuations of sums and products of Witt vectors of functions in the function field $K(C_k)$ of a curve $C_k$ over $k$ :

\medskip

\begin{lemma}
\label{sopro}
Let $f=(f_0,\dots,f_{l-1}),~g=(g_0,\dots,g_{l-1}) \in W_l(K(C_k))$, $P$ a point of $C_k(\bar{k})$, and $v_P$ the corresponding valuation. Set $a_i=v_P(f_i)~;~ b_i=v_P(g_i);~ 0\leq i\leq l-1.$ Then

i) if $h=(h_0,\dots,h_{l-1})=f+g$, we have :
$$v_P(h_i)\geq \min_{0\leq k \leq i}(p^{i-k}a_k,p^{i-k}b_k) ,\qquad 0\leq i\leq l-1.$$

ii) if $h=(h_0,\dots,h_{l-1})=fg$, we have :
$$v_P(h_0)=a_0+b_0,\qquad v_P(h_i)\geq \min_{0< j+k \leq i}(p^{i-j}a_j+p^{i-k}b_k),
\quad 1\leq i\leq l-1.$$
\end{lemma}

\section{Teichm\"uller subsets}

 We begin by recalling some facts about Greenberg's functor : the proofs can be found in \cite{gre}. To any $k$-scheme $Y$, we associate functorially a scheme $WY$ over $\Sp R$, with the same underlying topological space, and structure sheaf $\O_{WY}=W_l(\O_Y)$, defined for every $U$ open in $Y$ by :
$$\Gamma(U,\O_{W Y})=W_l(\Gamma(U,\O_Y)).$$
If $X$ is a $R$-scheme, the functor associating $\Ho_R(WY,X)$ to every $k$-scheme $Y$ is representable by a $k$-scheme $\FF X$ : we have a functorial isomorphism :
$$ \Ho_R(WY,X)=\Ho_k(Y,\FF X).$$
Now Greenberg's functor is the functor $X\mapsto \FF X$ from the category of $R$-schemes to the one of $k$-schemes. Note that we get a morphism $\lambda_X$ corresponding to the identity of $\FF X$ by the adjuction formula. Moreover, since we have $W\Sp k_d=\Sp R_d$ for any $d\geq 1$, this formula gives a bijection, we shall call {\it Greenberg's bijection} in the sequel :
$$X(R_d)=\FF X(k_d),$$
between the set of $R_d$-points of $X$ and the set of $k_d$-points of $\FF X$.

\medskip

We now describe the image of an affine $R$-scheme $X$ under the functor $\FF$, 
and explicit the morphism $\lambda_X$ in this case. If $X=\A_R^n$ is affine 
$n$-space over $R$, we have $\FF X=\A_k^{nl}$, the map $\lambda_X:=\lambda$ 
is the morphim corresponding to the morphism of $R$-algebras :
$$\begin{array}{ccccc}
\gamma & : & R[T_1,\dots,T_n] & \rightarrow & W_l\left(k[T_i^j]_{1\leq i \leq n,0\leq j\leq l-1}\right) \\
& & \alpha \in R & \mapsto & w_l(\alpha) \\
& & T_i & \mapsto & (T_i^0,\dots,T_i^{l-1}).\\
\end{array}$$

and Greenberg's bijection is just the map :

$$\begin{array}{ccccc}
 w_l\times \dots \times w_l&  : &\A^n_R(R_d)=R_d^n & \rightarrow & \A_k^{nl}(k_d)=k_d^{nl} \\
&& (\alpha_1\dots,\alpha_n)& \mapsto & (w_l(\alpha_1),\dots,w_l(\alpha_n)).\\
\end{array}$$
 
\medskip

Let us write the image $\gamma(f)$ of a polynomial $f \in R[T_i]$. We first fix some notations. Let $\delta$ be the morphism :
$$
\begin{array}{ccccc}
\delta & : & R[T_i] & \rightarrow & W_l(k[T_i^0])\\
 &  & T_i & \mapsto & (T_i^0,0,\dots,0)\\
\end{array}
$$
$\delta_{|R}=w_l$, and $U_i=(0,T_i^1,\dots,T_i^{l-1})$. Let $M=\{m_1,\dots,m_t\}$ denote a $t$-uple of integers greater than $1$ ; we set $|M|=m_1+\dots+m_t$, $\#M=t$, and $M!=m_1!\dots m_t!$. If $J=\{j_1,\dots,j_t\}$ is a $t$-uple of indexes in $\{1,\dots,n\}$, we set :
$$\frac{\partial^{|M|}f}{\partial T_J^M}=\frac{\partial^{|M|} f}{\partial T_{j_1}^{m_1}\dots \partial T_{j_t}^{m_t}}~;~U_J^M=U_{j_1}^{m_1}\dots U_{j_t}^{m_t}.$$

 Then we can write, for all  $f \in R[T_i]$:
\begin{eqnarray*}
 \gamma(f)
&
=
&
\delta(f)+\sum_{h=1}^{l-1}\sum_{M\atop{|M|=h}} \sum_{J\atop{\#J=\#M}}\delta\left(\frac{1}{M!}\frac{\partial^{h}f}{\partial T_J^M}\right)U_J^M.
\end{eqnarray*}
Note that $\frac{1}{M!}\frac{\partial^{h}f}{\partial T_J^M}$ is again a polynomial in $R[T_i]$. Developing, we obtain another description for $\gamma$ : if $\gamma(f)=(f_0,\dots,f_{l-1})$, we have :

$$f_0=\bar{f}(x_i^0), \qquad f_j=f^{(j)}(T_u^t)+
\sum_{i=1}^n (\frac{\partial \bar{f}}{\partial T_i}(T_i^0))^{p^{j}}T_i^{j},~1\leq j\leq l-1,$$

with $\bar{f}$ the reduction modulo $p$ of $f$, and $f^{(j)}$ a polynomial 
in $k[T_u^t]_{ 1\leq u \leq n,~ 0\leq t\leq j-1}$.

\medskip

Now let $U$ be an affine $R$-scheme $U=\Sp R[T_1,\dots,T_n]/I$, $I=(f_1,\dots,f_m)$. For every $1\leq i\leq m$, set $\gamma(f_i)=(f_{i0},\dots,f_{il-1})$ ; let $I'$ be the ideal in $k[T_u^t]$ generated by the $f_{ij}, 1\leq i \leq m, 0\leq j \leq l-1$, and $\pi_{I'}$ the surjection $k[T_u^t]\rightarrow k[T_u^t]/I'$. Clearly the morphism :
$$\xymatrix{
R[T_1,\dots,T_n] \ar[r]^{\gamma}& W_l(k[T_u^t]) \ar[r]^{W_l(\pi_{I'})} &
W_l(k[T_u^t]/I'),
}$$
vanishes on $I$ ; it induces a morphism $\gamma_U: R[T_1,\dots,T_n]/I \rightarrow  W_l(k[T_u^t]/I')$. Now the $k$-scheme $\Sp k[T_u^t]/I'$ is $\FF U$, and $\lambda_U$ is the morphism corresponding to $\gamma_U$.

\bigskip

In the sequel, we restrict our attention to the following situation : let $U=\Sp B$ be a smooth affine $R$-scheme, and $U_k=U\otimes k=\Sp B_0$ ; we denote by $\pi :U_k\rightarrow U$ the morphism induced by reduction modulo $p$ (note that since $U$ is a smooth $R$-scheme, $B$ is a flat $R$-algebra and $B_0=B\otimes k = B/pB$).  

\medskip

Let us construct a ``reduction morphism'' $\rho_U : \FF U\rightarrow U_k$. Projection of Witt vectors on their first coordinate, $p_0$, induces a morphism $\mu_{\FF U} : \FF U \rightarrow W\FF U$. Composing with $\lambda_U$, we get a morphism $\FF U \rightarrow U$, which factors to $\rho_U : \FF U \rightarrow U_k$. It corresponds to the morphism of $R$-algebras $p_0\circ \gamma_U$, which factors as :
$$B=R[T_1,\dots,T_n]/(f_1,\dots,f_m)\rightarrow B_0 \simeq k[T_1^0,\dots,T_n^0]/(\bar{f_1},\dots,\bar{f_m}) \hookrightarrow k[T_i^j]/I_{\FF U},$$
and the last inclusion is the morphism $\Gamma(\rho_U)$. Note that if $(a_1^0,\dots,a_n^{l-1})$ is a point of $\FF U$, $\rho_U$ sends it to $(a_1^0,\dots,a_n^0)$ : it acts as reduction modulo $p$ of the coordinates, seen as Witt vectors.

\medskip

We are interested in constructing sections of $\rho_U$. Note that such a section gives a lifting of the Frobenius morphism from $U_k$ to $U$ ; if $C$ is a smooth projective curve over $R$, with $C_k$ of genus at least $2$, this is impossible due to a result of Raynaud (\cite{ray} I.5.4). That is the reason why we restrict our attention to an affine open $U$. 

\medskip

We have $WU_k=\Sp W_l(B_0)$ : $U_k$ is defined by a nilpotent sheaf of ideals in $WU_k$. Since $U$ is smooth over $R$, the map $\Hom_{R}(WU_k,U) \rightarrow \Hom_{R} (U_k,U)$ induced by $\mu_{U_k}$ is surjective (\cite{ega} IV.17.1.1). Let $s$ be an antecedent of $\pi :U_k \rightarrow U$. The adjunction formula now gives $\sigma \in \Hom_k(U_k,\FF U)$ such that $\lambda_U\circ W\sigma=s$.

\medskip

{\bf Lemma 1.1.}{\it The morphism $\sigma$ is a section of $\rho_U$ ; it is a closed immersion $U_k \rightarrow \FF U$.}

\begin{proof}
Consider the following diagram : 

$$\xymatrix{
 \ar[d]_{\mu_{U_k}} U_k \ar[r]^{\sigma} & \ar[d]_{\mu_{\FF U}} \FF U \ar[r]^{\rho_U} & \ar[d]_{\pi} U_k \\
WU_k \ar[r]_{W\sigma} & W\FF U \ar[r]_{\lambda_U} & U \\
}$$
The left hand square is cartesian from the definition of the functor $W$, and the right hand one by the construction of $\rho_U$. From the adjunction formula, we have $\lambda_U\circ W\sigma=s$, and $\pi=s\circ \mu_{U_k}$ since $s$ is an antecedent of $\pi$ ; thus we get $\lambda_U\circ W\sigma \circ \mu_{U_k}=\pi$. Summing up, we obtain $\pi\circ \rho_U\circ \sigma=\pi$, and since $\pi$ is a closed immersion, we get the first claim of the lemma : $\rho_U\circ \sigma= \Id_{U_k}$. The second assertion comes immediately from the first one.
\end{proof}

Now we get a bijection $\sigma(U_k)(k)=U_k(k)$, and a bijection from $\sigma(U_k)(k) \subset \FF U(k)$ to a subset of $U(R)$.

\medskip

{\bf  Definition :}
\label{def}
{\it We call this subset the {\rm Teichm\"uller subset of $U(R)$} associated to the section $\sigma$, and denote it by  $\T_{\sigma}$.}

\medskip

{\bf Example :}
The simplest case is $U=\A_R^1={\rm Spec} R[T]$ ;
here we have $\FF U=\A_k^{l}=\Sp k[T^0,\dots,T^{l-1}]$, $\rho_U$ is the morphism $ k[T^0] \hookrightarrow  k[T^0,\dots,T^{l-1}]$, and $\sigma :  k[T^0,\dots,T^{l-1}] \rightarrow k[T]$,
$T^0 \mapsto T$, $T^i\mapsto 0, \quad 1\leq i\leq l-1$ is a section of $\rho_U$ ; the associated Teichm\"uller subset (via the identification $\A_R^1(R)=R$), is just the Teichm\"uller $\T$ of $R$.

\medskip

Now we come to the main property of Teichm\"uller subsets ; they allow us to rewrite exponential sums :

\medskip

{\bf Theorem 1.2.}
\label{sig}
{\it Let $f\in B$ be a regular function over $U$, and $\psi$ an additive character of
$GR(p^l,m)$. We have the equality :

$$\sum_{\Pi\in \T_{\sigma}} \psi\left(f(\Pi)\right) =
 \sum_{P\in U_k(k)} \psi\left(\Gamma(s)(f)(P)\right), $$
where $\Gamma(s)(f)(P)=(f_0(P),\dots,f_{l-1}(P))$ if $\Gamma(s)(f)=(f_0,\dots,f_{l-1})$}

\begin{proof}
Let $P \in U_k(k)$, $\tilde{P}=\sigma P \in \FF U(k)$, and $\Pi \in
\T_\sigma$ the point corresponding to $\tilde{P}$ by Greenberg's bijection. Consider the following diagram (where $P,\tilde{P},\Pi$ stand for the morphisms associated to the points $P,\tilde{P},\Pi$) :  
$$\xymatrix{
 & & \ar[dll]_{\Gamma(s)} B \ar[d]_{\Gamma(\lambda_U)} \ar[ddrr]^{\Gamma(\Pi)} & &\\
W_l(B_0) \ar[drr]_{W_l(\Gamma(P))} & & \ar[ll]_{W_l(\Gamma(\sigma))} \ar[d]^{W_l(\Gamma(\tilde{P}))}
W_l(\Gamma(\FF U,\O_{\FF U})) & &\\ 
& & W_l(k) \ar[rr]^{w^{-1}} & & R
}$$
Since $\tilde{P}=\sigma P$, we have the equality of morphisms $\sigma \circ P= \tilde{P}$, that is $\Gamma(\tilde{P})= \Gamma(P)\circ \Gamma(\sigma)$ ; moreover, since $\tilde{P}$ and $\Pi$ correspond by Greenberg's bijection, the right-hand side of the diagram commutes. Thus the whole diagram commutes and we get :
$$\forall f\in B,~\Gamma(s)(f)(P)=w(f(\Pi)).$$
Finally the two summation sets are in one to one correspondance,
and the theorem is established.
\end{proof}

\section{The degrees of coordinate functions}

In this section, we give a precise description of certain sections $\sigma$, when $U$ is an affine open of a smooth curve $C\subset \P_R^n$.

\medskip

Let $C$ be a smooth, irreducible, geometrically connected curve over $R$, equiped with an embedding $C\subset \P_R^n=\Proj R[X_0,\dots,X_n]$, and $U$ the affine open $C\cap \{X_0\neq 0\}$ of $C$. Recall that if $C$ is defined in $\P_R^n$ by homogeneous polynomials $F_1,\dots,F_m$ with $\deg F_i=d_i$, $U$ is defined in $\A_R^n=\Sp R[x_1,\dots,x_n]$ by $I=(f_1,\dots,f_m)$, $f_i=F_i/X_0^{d_i}$. We set $B=R[x_1,\dots,x_n]/I$, and $B_0=B\otimes k$ the ring of regular functions on $U_k=U\otimes k$.

\medskip

Recall that the morphism : $\Gamma(\rho_U):B_0 \rightarrow k[x_i^j]_{1\leq i\leq n,~0\leq j\leq l-1}/I_{\FF U}$ is the injection $x_i\mapsto x_i^0$. Thus a morphism $\phi:k[x_i^j]_{1\leq i\leq n,~0\leq j\leq l-1}\rightarrow B_0$ factors to a section $\sigma$ of $\rho_U$ if and only if $\phi(x_i^0)=x_i$ and it vanishes on $I_{\FF U}$. Since the generators of $I_{\FF U}$ are the components of the Witt vectors $\gamma(f_i)$, this last condition can be rewritten : 
$$\phi(\bar{f_t}(x_1^0,\dots,x_n^0))=0 \quad  1\leq t \leq m~; $$
$$\phi \left(\sum_{i=1}^n (\frac{\partial \bar{f_t}}{\partial x^0_i})^{p^j}x_i^{j}
+f_t^{(j)}\right)=0  \quad 1\leq t \leq m, \quad 1\leq j\leq l-1.$$

in $B_0$. Now since $\phi(x_i^0)=x_i$, the first conditions are automatically satisfied, and, if we set $\phi(x_i^j)=x_{ij}$, we can rewrite the second ones in $B_0$ as :
$$\sum_{i=1}^n (\frac{\partial \bar{f_t}}{\partial x_i})^{p^j}x_{ij}+f_t^{(j)}(x_i,x_{i1},\dots,x_{i,j-1})=0$$

\medskip

{\bf Definition :} If $\sigma$ is the corresponding section of $\rho_U$, we call the $(x_{ij})$  the {\it coordinate functions} of the Teichm\"uller subset $\T_\sigma$. 

\medskip

We now give a bound for the pole orders of these functions.

\medskip

{\bf Theorem 2.1. }
\label{deg}
{\it Let $C$, $U$ be as above, assume that $C_k$ has genus $g$, and denote by $D_0$ the intersection divisor of $C_k$ with the hyperplane $\{ X_0=0\} \subset \P_k^n$. Then if $D$ is an effective divisor on $C_k$ of degree $2g-2+{\lceil\frac{2g-1}{p}\rceil}$ whose support is contained in that of $D_0$, we can find a section $\sigma : U_k\rightarrow \FF U$ of $\rho_U$ with coordinate functions $(x_{ij})$ such that, in $W_l(B_0)$ :
$$(0,x_{i,1},\dots,x_{i,l-1})=\sum_{d=1}^{l-1}p^{d-1}V\left(X_i^{(d)}-(Fx_i)X_0^{(d)}\right),~1\leq
 i\leq n,$$
with $X_i^{(d)}=(x_{i,d}^{(d)},\dots,x_{i,l-1}^{(d)},0,\dots,0)$ satisfying $x_{i,j}^{(d)}\in L\left(p^{j-d+1} d(D+D_0)\right)$}.

\medskip

We need a lemma :

\medskip

{\bf Lemma 2.2. }
\label{tec}
{\it Let $d,e\in \N^*$, and $E=dD+(d-1)D_0$. If $R\in \Gamma(U_k,F^{*e}\n_{\ma{P}_k^n/C_k})$ is a global section of the sheaf $F^{*e}\left(\n_{\ma{P}_k^n/C_k}\otimes \L(E )\right)$, there exists a global section $\eta$ of $F^{*e}\left(\T_{\ma{P}^n}\otimes \L(E)\right)$ whose image is $R$ by the map :
$$\Gamma\left( C_k,F^{*e}\left(\T_{\ma{P}^n}\otimes \L(E)\right)\right) \rightarrow \Gamma\left(C_k,F^{*e}\left(\n_{\ma{P}_k^n/C_k}\otimes \L(E)\right)\right),$$
and whose restriction :
$$(\eta_1,\dots,\eta_n) \in \Gamma\left(U_k, F^{*e}\left(\T_{\ma{P}^n}\otimes \L(E)\right)\right)= B_0^n$$
can be written :
$$\eta_i=E_i-x_i^{p^e}E_0,\quad E_i \in L(p^e(E+D_0)),~0\leq i\leq n.$$}

\begin{proof}
Consider the normal/tangent exact sequence of $C_k$ in $\P_k^n$ :
$$\xymatrix{
0\ar[r] & \T_{C_k} \ar[r] & \T_{\ma{P}^n_k}\otimes \O_{C_k} \ar[r] & \n_{\ma{P}_k^n/C_k} \ar[r] & 0
}$$

If we tensor it by $\L(E)$, then take its pull-back by the $e$-th power of the absolute Frobenius morphism, the sequence remains exact.
Now by the assumption on $D$, the first sheaf is invertible of degree
greater than $2g-1$, and its first cohomology group vanishes from
Riemann-Roch theorem. Hence the sequence remains exact when we apply the global
section functor. From this we deduce the existence of $\eta$.

\medskip

Consider now the pull-back of Euler's exact sequence to $C_k$ :
$$\xymatrix{
0\ar[r] & \O_{C_k} \ar[r] & \O_{C_k}(1)^{n+1} \ar[r] & \T_{\ma{P}^n_k}\otimes \O_{C_k} \ar[r] & 0
}$$ 
Remark that $\O_{C_k}(1)\simeq \L(D_0)$ and apply the same operations as above.
Again, the sequence remains exact when we take global sections, and we obtain a preimage for $\eta$ :
$$E=(E_0,\dots,E_n) \in \Gamma\left(F^{*e}\L(D_0+E)\right)^{n+1}=L(p^e(D_0+E))^{n+1}.$$
For the last assertion, we just have to apply the (exact) functor $\Gamma(U_k,.)$ to the last exact sequence ; we get :
$$\Gamma\left(U_k, F^{*e}\L(D_0+E)^{n+1}\right)=(B_0^{*e})^{n+1}$$
$$\Gamma\left(U_k, F^{*e}(\T_{\ma{P}^n_k}\otimes \L(E))\right)=(B_0^{*e})^{n}$$
(where $B_0^{*e}$ means the $B_0$-module $B_0$ acting on itself via the $e$-th power of the Frobenius morphism), and the morphism from the first of these $B_0$-modules to the second is :
$$(b_0,\dots,b_n)\mapsto (b_1-x_1^{p^e}b_0,\dots,b_n-x_n^{p^e}b_0),$$
This completes the proof of lemma 2.2.
\end{proof}

\begin{proof} (of theorem 2.1.)

Let us begin by the case $l=2$. Consider the morphism :
$$\begin{array}{ccccc}
\phi_1 & : & R[x_1,\dots,x_n] & \rightarrow & W_2(B_0) \\
& & x_i &\mapsto & (x_i,0)
\end{array}$$
Since $\phi_1$ induces reduction modulo $p$ on the first coordinate, we see that $\phi_{1|I}$ induces a $B$-module homomorphism $S_1:I/pI \rightarrow V(B_0)\subset W_2(B_0)$ such that $(I/pI)^2 \subset \Ker S_1$. Note that since $I$ is the kernel of $R[x_1,\dots,x_n]\rightarrow B$, it is a flat $R$-module ; thus $I/pI=I_0$, and we get a morphism :
$$R_1 :I_0/(I_0)^2 \rightarrow V(B_0)\simeq B_0^*,$$
where $B_0^*$ is the $B_0$-module described in the proof of the lemma. It is therefore natural to consider $R_1$ as a section of $\Gamma(U_k, F^*\n_{\ma{P}_k^n/C_k})$.
Now if $\tilde{f_i}$ is the image of $\bar{f_i}$ in $I_0/(I_0)^2$, then  $f_t^{(1)}(x_i)=R_1(\tilde{f_t})$, and solving the equations (in $B_0$) :
$$\sum_{j=1}^n (\frac{\partial \bar{f_t}}{\partial x_j})^{p}x_{1j}+f_t^{(1)}(x_j)=0
\quad 1\leq t \leq m$$
boils down to finding a preimage of $-R_1$ in $\Gamma(U_k,F^*(\T_{\ma{P}_k^n}\otimes
\O_{C_k}))$ by the map :
$$\Gamma(U_k,F^*(\T_{\ma{P}_k^n}\otimes
\O_{C_k})) \rightarrow \Gamma(U_k,F^*\n_{\P_k^n/C_k}),$$
which is the morphism of $B_0^*$-modules induced by the
transpose of the Jacobian matrix of $U_k$ in $\A_k^n$. Let us show the :

\begin{claim}
$R_1$ is a global section of $F^*\n_{\ma{P}_k^n/C_k}$.
\end{claim}

It is sufficient to show that $R_1 \in (F^*\n_{\ma{P}_k^n/C_k})_P$ for all $P$ not in $U_k$ since $R_1 \in \Gamma(U_k, F^*\n_{\ma{P}_k^n/C_k})$. Let $P$ be such a point, and assume it lies away from the hyperplane $\{X_n = 0\}$. By the Jacobian criterion of smoothness, the $\O_{C_k,P}$-module $(I_0/I_0^2)_P$ is generated by $n-1$ functions :
$$\frac{\overline{F_1}}{X_n^{d_1}}=\frac{\tilde{f_1}}{x_n^{d_1}}, \dots,\frac{\overline{F_{n-1}}}{X_n^{d_{n-1}}}=\frac{\widetilde{f_{n-1}}}{x_n^{d_{n-1}}},~d_i=\deg f_i.$$
Their images by $R_1$ are :
$$\frac{R_1\tilde{f_1}}{x_n^{pd_1}}, \dots,\frac{R_1\widetilde{f_{n-1}}}{x_n^{pd_{n-1}}}~;$$
by Witt vectors calculations we get $R_1\tilde{f_r}=\frac{1}{p}(f_r(x_i)^p-f_r^p(x_i^p))$, and $\deg R_1\tilde{f_r}\leq pd_r$. Finally the $\frac{R_1\tilde{f_r}}{x_n^{pd_r}}$ are in $\O_{C_k,P}$, and $R_1 \in (F^*\n_{\ma{P}_k^n/C_k})_P$ : this proves the claim.

\medskip

Consequently, since $D$ is an effective divisor with support contained in
$C_k\backslash U_k$, we get that $R_1$ is also a global section of the
sheaf $F^*(\n_{\ma{P}_k^n/C_k}\otimes\L(D))$. Now we can apply lemma 2.2
to $-R_1$ with $d=1, e=1$ : we obtain functions $x_{i,1}^{(1)},~0\leq i\leq n$ that are global sections of $\L(p(D+D_0))$, and theorem 2.1 is proved for $l=2$.

\bigskip

We finish the proof by induction on $l$ : suppose we have shown the result for $l-1$, and consider the morphism :
$$\begin{array}{ccccc}
\phi_{l-1} & : & R[x_1,\dots,x_n] & \rightarrow & W_l(B_0) \\
& & x_i & \mapsto & (x_i,0,\dots,0)+\sum_{k=1}^{l-2}p^{k-1}V(Y_i^{(k)}-(Fx_i)Y_0^{(k)})
\end{array}$$
where $Y_i^{(k)}=(x_{i,k}^{(k)},\dots,x_{i,l-2}^{(k)},0,\dots,0)$ with the $x_{i,j}^{(k)}$ as in theorem 2.1 for $l-1$. As above, the morphism 
$\phi_{l-1|I}$ factors to a morphism of $B_0$-modules :
$$R_{l-1}:I_0/I_0^2 \rightarrow V^{l-1}(B_0)\simeq B_0^{*l-1},$$
 $R_{l-1}(\tilde{f_t})=f_t^{(l-1)}(x_i,x_{ij})$ : we look for a preimage. For this purpose we need to rewrite $f_t^{(l-1)}$  more precisely ; from the description of $\gamma$, we get :
$$\phi_{l-1}(f)=\delta(f)+\sum_{h=1}^{l-1}\sum_{|M|=h\atop{\#J=\#M}}
\delta\left(\frac{1}{M!}\frac{\partial^h f}{\partial x_J^M}\right)
\prod_{j_i\in JM}\left(\sum_{k=1}^{l-2}p^{k-1}V(Y_{j_i}^{(k)}-Fx_{j_i}Y_0^{(k)})\right),$$
where $JM$ is the set of the indexes $j_u$ in $J$, each one counted with the 
multiplicity $m_u$ (note that $JM$ contains $h$ elements). Developing the last product yields the following general term : 
$$p^{k_1+\dots+k_h-h}\prod_{j_i\in JM} V(Y_{j_i}^{(k_i)}-(Fx_{j_i})Y_0^{(k_i)})
=p^{k_1+\dots+k_h-1}V\left(\prod_{j_i\in JM} (Y_{j_i}^{(k_i)}-(Fx_{j_i})Y_0^{(k_i)})\right)$$
Therefore we can write :
$$
\phi_{l-1}(f)=R_{l-1}^{(1)}(f)+\sum_{d=2}^{l-1} p^{d-1} R_{l-1}^{(d)}(f),$$
$$R_{l-1}^{(1)}(f)=\delta(f) +\sum_{i=1}^n \delta\left(\frac{\partial f}{\partial x_i}\right) V(Y_i^{(1)}-(Fx_i)Y_0^{(1)}),~\dots$$
$$R_{l-1}^{(d)}(f)= \sum_{h=1}^{d}\sum_{k_1+\dots+k_h=d} 
V\left(\sum_{|M|=h\atop{\#J=\#M}} 
F\delta\left(\frac{1}{M!}\frac{\partial^h f}{\partial x_J^M}\right)
\prod_{j_i\in JM} (Y_{j_i}^{(k_i)}-(Fx_{j_i})Y_0^{(k_i)})\right)$$

Now by the induction hypothesis, for $f\in I$ we have $R_{l-1}^{(d)}(f) \in V^{l-d}(B_0)\subset W_{l-d+1}(B_0)$ and as above,  $R_{l-1}^{(d)}$ factors to a morphism from $I_0/I_0^2$ to $V^{l-d}(B_0)\simeq B_0^{*l-d}$ ; thus we write :
$$R_{l-1}^{(d)} \in \Gamma\left(U_k, F^{*l-d}\n_{\ma{P}_k^n/C_k}\right).$$

\begin{claim}
$R_{l-1}^{(d)}$ is a global section of the sheaf $ F^{*l-d}(\n_{\ma{P}_k^n/C_k}\otimes \L(dD+(d-1)D_0))$.
\end{claim}

As in the proof of claim 1, let $P \in C_k\backslash U_k(\bar{k})$ lying away from
the hyperplane $\{X_n=0\}$, and $f_r/x_n^{d_r},~1\leq r\leq n-1$, 
a set of generators for the free $\O_{C_k,P}$-module $(I_0/I_0^2)_P$ ; 
we can write $R_{l-1}^{(d)}(f_r/x_n^{d_r})$ as :  

$$\frac{1}{x_n^{p^{l-d}d_r}} \sum_{h=1}^{d}\sum_{k_1+\dots+k_h=d} 
V\left(\sum_{|M|=h\atop{\#J=\#M}} F\delta\left(\frac{1}{M!}
\frac{\partial^h f_r}{\partial x_J^M}\right)\prod_{j_i\in JM} (Y_{j_i}^{(k_i)}
-(Fx_{j_i})Y_0^{(k_i)})\right).$$

When we develop the last product in the numerator, the terms corresponding 
to the systematic choice of $(Fx_{j_i})Y_0^{(k_i)}$ give :
$$Y_0^{(k_1)}\dots Y_0^{(k_h)}F\delta\left(\sum_{|M|=h\atop{\#J=\#M}}
\frac{1}{M!}\frac{\partial^h f_r}{\partial x_J^M} x_J^M\right),$$
with $1\leq h \leq l-1$ and $k_1,\dots,k_h$ as above. 
But the function in the parentheses is the image in $B_0$ of a polynomial of 
degree less than $d_r-1$ (its homogeneous part of degree $d_r$ is a 
multiple of that of $f_r\in I$ by Euler's theorem). Denoting by $v$ 
({\it resp.} $v_0$) the multiplicity of $D$ ({\it resp.} $D_0$) at $P$, and applying lemma 0.1,
the $u$-th component of this term has valuation at $P$ at least :
$$-p^{u+1}(k_1+\dots+k_h)(v+v_0)-p^{u+1}(d_r-1)v_0=-p^{u+1}d_rv_0-p^{u+1}(dv+(d-1)v_0).$$
The $u$-th components of the other terms of the numerator (where we choose 
at least once $Y_{j_i}^{(k_i)}$ in the development) have valuation at $P$ at least :
$$-p^{u+1}\left((d_r-h)v_0+\sum_{i=1}^h (k_i(v+v_0) + \epsilon_iv_0)\right)
\geq -p^{u+1}d_rv_0-p^{u+1}(dv+(d-1)v_0)$$
where $\epsilon_i\in \{0,1\}$ and at least one of the $\epsilon_i$ is zero 
(note that $\epsilon_i=1$ corresponds to the choice of $(Fx_{j_i})Y_0^{(k_i)}$ 
rather than $Y_{j_i}^{(k_i)}$ in the
product). Since the sum of these vectors is in $V^{l-d-1}(B_0)$, we write it 
$V^{l-d-1}(F_r)$, and we get $v_P(F_r)\geq -p^{l-d}(d_r v_0+dv+(d-1)v_0)$. Finally, 
since the denominator $x_n^{p^{l-d}d_r}$ has valuation exactly $-p^{l-d}d_rv_0$, this ends the proof of the claim.

\medskip

Applying lemma 2.2 with $e=l-d$, we deduce the existence of $x_{i,l-1}^{(d)}\in \L(p^{l-d}d(D+D_0)),~0\leq i\leq n$ such that $-R_{l-1}^{(d)}$ is the image of :
$$(x_{i,l-1}^{(d)}-x_i^{p^{l-d}}x_{0,l-1}^{(d)})_{1\leq i \leq n} \in \Gamma\left(U_k,F^{*l-d}(\T_{\ma{P}_k^n}\otimes\L(dD+(d-1)D_0))\right).$$
and setting : 
$$ X_i^{(d)}=Y_i^{(d)}+V^{l-1-d}(x_{i,l-1}^{(d)},0,\dots,0)~;~X_i^{(l-1)}=(x_{i,l-1}^{(l-1)},0,\dots,0)$$
gives a Witt vector satisfying the hypotheses of theorem 2.1. Finally, the morphism defined by :
$$\phi_{l-1}(x_i)=(x_i,0,\dots,0)+\sum_{d=1}^{l-1}p^{d-1}V(X_i^{(d)}-(Fx_i)X_0^{(d)}),$$
vanishes on $I$ from the preceding construction, and therefore defines a section of $\rho_U$. This completes the proof of theorem 2.1.
\end{proof}

The following is an easy consequence of the theorem, and will be useful in the proof of Proposition 4.3.

\smallskip

{\bf Corollary 2.3.}
{\it The function $x_{ij}$ is in $L(p^j(jD+(j+1)D_0))$.}

\medskip

{\bf Remark :}
The degrees of the coordinate functions of the elliptic Teichm\"uller
lift of an ordinary elliptic curve over a perfect field of characteristic
$p\neq 2,3$ have been obtained by Finotti in
\cite{fin} ; the exact degrees are obtained by calculations from
the equation $y^2=x^3+ax+b$, and using the group structure.

\section{Exponential sums}

The aim of this section is to recall bounds for exponential sums over a curve $C_k$, associated to an additive character of a Galois ring. We first describe Artin-Schreier-Witt theory, i.e. the theory of abelian $p$-extensions of a field of characteristic $p$, then we focus on algebraic function fields of one variable before giving the bound. 
\medskip

Let $K$ be a field of characteristic $p$, fix $K_s$ a separable closure of $K$, and let $G=\Gal(K_s/K)$ be the absolute Galois group of $K$. We define the morphism $\wp=F-\Id$ on $W_l(K_s)$. It is a $G$-morphism, surjective on $W_l(K_s)$, and we know from Galois theory that there is a one to one inclusion preserving correspondance between the collection of subgroups of $W_l(K)/\wp W_l(K)$ and the collection of abelian extensions of $K$ of exponent a divisor of $p^l$.

\medskip

If $A$ is a subgroup of $W_l(K)/\wp W_l(K)$, let us denote by $K(\wp^{-1}A)$ the corresponding extension. In the sequel, we restrict our attention to the extensions $K(\wp^{-1}<f>):=K(\wp^{-1}f)$ where $K$ is an algebraic function field of one variable, the function field of a curve $C_k$, with field of constants $k=\F_{p^m}$, and $f$ an element of $W_l(K)$. 

\medskip

{\bf Definition.}
{\it a) Let $d$ be an integer, $P$ be a place of $K$ (a point of $C_k(\bar{k})$), and $v$ the corresponding valuation. We define a subset $V_{d,P}\subset W_l(K)$ by 
$$V_{d,P}=\{ (a_0,\dots,a_{l-1})\in W_l(K)/~ \forall i\in\{0,\dots,l-1\} ,
~ p^{l-1-i}v(a_i)\geq d\}.$$

b) If $f\in W_l(K)$, we define the {\rm valuation of $f$ at $P$} by $V_P(f)=\max\{d\in \Z, f\in V_{d,P}\}$, and the {\rm reduced pole order of $f$ at $P$} by :
 $$\rp_P(f)=\min\{ d\geq -1,f\in V_{-d,P}+\wp W_l(K)\}.$$}

If $a,b$ are in $W_l(K)$, we have $\rp_P(a+b) \leq \max(\rp_P(a),\rp_P(b))$, with equality if $\rp_P(a)\neq \rp_P(b)$ ; moreover, if $0\leq k \leq l-1$, we have : $p^k\rp_P(V^k a) \leq  \rp_P(a)$.

\medskip

At a fixed place $P$, we can always reduce a Witt vector modulo $\wp W_l(K)$ in such a way that none of its components has its valuation divisible by $p$ : if $f$ is in $W_l(K)$, there exists $f'=(f'_0,\dots,f_{l-1}')\in W_l(K)$ such that :
$$i) f\equiv f' [\wp W_l(K)]~ ;\quad
ii) v_P(f_i')\geq 0 \text{ or } (v_P(f_i'),p)=1,~0\leq i\leq l-1.$$
This is the equivalent of Artin's algorithm to compute reduced pole orders in the case $l=1$ (cf \cite{lac} Proposition 1.8) ; if $\rp_P(f)>0$, we have :
$$\rp_P(f)=\max_{0\leq i\leq l-1}\{-p^{l-1-i}v_P(f_i')\}.$$

As in the Artin-Schreier case, the extension  $K(\wp^{-1}f)$, $f=(f_0,\dots,f_{l-1})$ is of degree $p^l$, and its constant field is $k$, if and only if $f_0 \notin k+\wp_0K$, where $\wp_0$ is the map $g\mapsto g^p-g$ from $K$ to $K$. We say that $f$ is {\it nondegenerate}. In this case we have the following estimate for the genus of $L=K(\wp^{-1}f)$ :
$$2(g_L-1) = 2p^l(g_K-1)+\sum_{n=1}^{p^l-1} \sum_{P \in \ma{P}(nf)} (\rp_P(nf)+1)\deg P,$$
where $\ma{P}(f)=\{ P\in C_k(\bar{k}),\quad \rp_P(f) >0 \}$.

\medskip

We come to exponential sums. First we define a local symbol attached to the extension
$K(\wp^{-1}f)/K$. We describe it explicitely from its ghost components (cf. \cite{wi}) : if $g\in K_P$, the completion of $K$ at $P$, $(f,g)_P$ is the element of $W_l(k(P))$, $k(P)$ the residue  field of $C_k$ at $P$, whose ghost components are :
 $$\Res_P(f_0\frac{dg}{g}),\dots,\Res_P((f_0^{p^{l-1}}+\dots+p^{l-1}f_{l-1})\frac{dg}{g}).$$

From an additive character of $W_l(k)$, $\psi$, we can define an idele class character : if $g=(g_P)_{P\in C_k(\bar{k})}$is an idele, we set :
$$\chi(g)=\prod_{P\in C_k(\bar{k})} \psi\left(\Tr_{W_l(k(P))/W_l(k)}(f,g_P)_P\right).$$
Its conductor is given by :
$$ D_{\chi} = \sum_{P\in \ma{P}(f)} (\rp_P(f)+1) P,$$
and we associate to this character a $L$-function :
$$L_{K,\chi}(T)=\prod_{P\in C_k(\bar{k})\backslash \ma{P}(f)} \frac{1}
{1-\psi\left(\Tr_{W_l(k(P))/W_l(k)}(f,\pi_P)_P\right)T^{\deg P}},$$
with $\pi_P$ a local parameter at $P$. This $L$-function is actually a polynomial of degree $\deg D_{\chi}+2g_K-2$, and the Riemann-Weil hypothesis for function fields allows us to give a bound for the following exponential sums, which are the coefficients of its development in power series :

\medskip

{\bf Theorem 3.1. }
{\it Let $f\in W_l(K)$ be a nondegenerate vector. We have :
$$ \left|\sum_{P\in C_k(k_d) \backslash \ma{P}(f)} \psi(\Tr_{W_l(k_d)/W_l(k)}f(P))\right|
\leq \left(  2(g_K-1)+\sum_{P \in \ma{P}(f)} (\rp_P(f)+1)\deg P\right) p^{\frac{md}{2}}.$$}

\section{Estimation of reduced pole orders}

In the sequel, $U=\Sp B$ is an affine open of a smooth projective $R$-curve $C$, $\T_\sigma$ is the Teichm\"uller subset associated to a section $\sigma : U_k\rightarrow \FF U$ satisfying theorem 2.1, and $\Gamma(s):B\rightarrow W_l(B_0)$ is the morphism of $R$-algebras associated to the morphism $s:WU_k\rightarrow U$ corresponding to $\sigma$.
From theorems 1.2 and 3.1, in order to give a bound for the exponential sum associated to a function $f\in B$ over the Teichm\"uller subset $\T_\sigma$, we have to estimate the reduced pole orders of the vector $\Gamma(s)(f)$ from the pole orders of $f$. We need to define what is the pole order of a function on $C$, a curve over $R$, at a point $\Pi$ of $C(R)$ ; this will be done from the expansion of $f$ in powers of local parameters, which are defined in section 4.1. Then we construct the morphism $S$ extending $\Gamma(s)$ to the (constant) sheaf of total quotient rings of $\O_C$ in section 4.2 ; finally we estimate reduced pole orders of the Witt vectors $S(t^i)$, where $t$ is a local parameter at a pole of $f$, and we obtain the reduced pole order of $\Gamma(s)(f)$ from its expansion in powers of $t$ in section 4.3.

\bigskip

\subsection{Local parameters}

Let $C$ be an irreducible, geometrically connected smooth of relative dimension $1$ scheme $C$ over $\Sp(R)$, equiped with an embedding $C\subset \P_R^n$.

Let $C_k:=C\times_{\Spec R} \Spec k$ be the ``reduction modulo $p$'' of $C$ ; it is a smooth irreducible projective curve over $\Spec k$. Note that, for any $d\geq 1$, points of $C(k_d)$ are in one to one correspondance with points of $C_k(k_d)$. We identify these two sets in the sequel.

\medskip
 
We denote by $K_C$ the sheaf of total quotient rings of $\O_C$ ; recall from \cite{wal} that it is a constant sheaf, whose associated ring we also denote by $K_C$.

\medskip

Let us denote by $C(R_d)$ the set of $R$-morphisms $\Pi :\Spec R_d\rightarrow C$ ; we call such a morphism an $R_d$-point of $C$. The morphism $\Sp k_d \rightarrow \Sp R_d$ induced by the reduction modulo $p$ associates to every $R_d$-point $\Pi$ an element $P \in C(k_d)$ ; we get a map $r_C:C(R_d)\rightarrow C(k_d)=C_k(k_d)$. We say : ``$\Pi$ is an $R_d$-point above $P$''.

\medskip

The ideal sheaf associated to an $R_d$-point $\Pi$ of $C$ (considered as a subscheme of $C$) is locally principal \cite{wal} thus there exists an element $t \in \O_{C,P}$ such that the kernel of the morphism : $\Pi :  \O_{C,P} \rightarrow R_d$ is the principal ideal $(t)\subset \O_{C,P}$. We call such an element a {\it local parameter} for $C$ at $\Pi$. If $m_P$ denotes the maximal ideal of $\O_{C,P}$, we have the following properties for a local parameter $t$ (cf.  \cite{wal}) :

i) $t$ is not a zero divisor in $\O_{C,P}$ ;

ii) $t \in m_P\backslash m_P^2$ ;

iii) $\O_{C,P}/(t^n)$ is a free $R_d$-module with basis $1,t,t^2,\dots,t^{n-1}$ ;

iv) the quotient ring of $\O_{C,P}$, $K_C$, is an $\O_{C,P}$-module generated by $1,t^{-1},\dots$.

\medskip

In the sequel we need to know local parameters for $R_d$-points explicitely : this is the aim of the next proposition.

\medskip

{\bf Proposition 4.1. }
\label{rpts}
{\it Let $C\subset \P^n_R=\Proj R[X_0,\dots,X_n]$ be a smooth projective curve over $R$, and $P$ be a point of $C(k_d)$, not lying in the hyperplane $\{X_n=0\}$. Then there exists an $R_d$-point $\Pi$ above $P$ whose local parameter is the image in $K_C$ of the function :
$$t=\pi_d(\frac{X_i}{X_n}),$$
where $i\in \{0,\dots,n-1\}$, and $\pi_d$ is an irreducible polynomial of degree $d'$ dividing $d$, with all its roots in $\T_{d'}$}.

\begin{proof}
This is a local question, so we restrict our attention to the open $U_0=\Sp R[x_0,\dots,x_{n-1}]\subset \P_R^n$ (with $x_i=X_i/X_n$) containing $P$. Let $C\cap U_0=U=\Sp R[x_0,\dots,x_{n-1}]/I_0$, where $I_0=(f_1,\dots,f_m)$ is the ideal defining $U\subset U_0$.
By the Jacobian criterion of smoothness applied to $U$ at $P$, we get $n-1$ functions among the $f_i,~1\leq i\leq m$ and $n-1$ indexes in $\{0,\dots,n-1\}$. Let $i$ be the last element of this set, and $\pi_d$ the polynomial with roots in $\T_d$, whose reduction modulo $p$ is the minimal polynomial of the $i$-th coordinate of any one of the $d$ points $P_d$ in $\P^n_{R_d}(k_d)$ corresponding to $P$. Consider the closed subscheme $S\subset U_0$ defined by the ideal $(f_1,\dots,f_{n-1},\pi_d(x_i))$. Thus $P$ is an element of $S(k_d)$, and the Jacobian criterion ensures us that $S$ is smooth over $\Sp R$ at $P$, of relative dimension $0$. We deduce that $\O_{S,P}$ is a local ring of dimension $0$, and a flat $R$-module (hence free since $R$ is a local ring), with maximal ideal $(p)$ and residue field $k_d$. Therefore we must have $\O_{S,P}\simeq \O_{C,P}/(\pi_d(x_{i}))_{P} \simeq R_d$, and the proposition follows.

\end{proof}

There is an other description of $R_d$-points, in term of Cartier divisors (cf \cite{wal}) : let $\Pi$ be an $R_d$-point of $C$, above $P\in C(k_d)$, and $t_d$ a local parameter for $C$ at $\Pi$. We associate to $\Pi$ the Cartier divisor $(\Pi)$ defined by : $(\Pi)=\left( (U,t_d), (V,1) \right)$,
where $V$ is the open $C\backslash \{P\}$, and $U$ an open containing $P$. We denote by $\O_C(\Pi)$ the invertible sheaf corresponding to this Cartier divisor.

\bigskip

\subsection{The morphism $S$}
We will work with expansions of $f$ at its poles, in terms of the local parameters obtained in proposition 1.1. Since these parameters are not in $B$, we first show that the morphism $\Gamma(s)$ extends to a morphism $S:K_C \rightarrow W_l(K(C_k))$ with $K_C$ as in section 1, and $K(C_k)$ the function field of the curve $C_k$.

\medskip

Note from the definition of $K_C$ that we have $K_C=i(U)^{-1}\Gamma(U,\O_C)=i(B)^{-1}B$, where $i(B)$ is the set of
elements of $B$ which are not zero divisors in any of the local
rings $B_{\mathfrak p}$. These are just the elements of
$B\backslash pB$ (since $\Sp B \rightarrow \Sp R$ is smooth, the
zero divisors must be nilpotent), and we get $B_{(p)}=K_C$. If $g\in B\backslash pB$, $\Gamma(s)(g)$ is an element of $W_l(B_0)$ with nonzero first coordinate, thus invertible in $W_l(K(C_k))$. We now define $S$ on $K_C$ by 
$$S(fg^{-1})=\Gamma(s)(f)\left(\Gamma(s)(g)\right)^{-1}, \quad f\in B,\quad g\in B\backslash pB.$$

Moreover, if $f,g$ are the images of $f'\in R[x_1,\dots,x_n]$, $g'\in R[x_1,\dots,x_n]\backslash(p)$, we can extend $\delta$ to $R[x_1,\dots,x_n]_{(p)}$, and we can write again, in $W_l(K(C_k))$ :

$$S(\frac{f}{g})=
\delta(f'/g')+\sum_{h=1}^{l-1}\sum_{|M|=h,\atop{\#J=\#M}}
\delta\left(\frac{1}{M!}\frac{\partial^h (f'/g')}{\partial x_J^M}\right)U_J^M.$$

\bigskip

\subsection{Reduced pole orders}

If $f\in \Gamma(U,\O_C)$ is a regular function on $U$, we have to estimate reduced pole orders of the
function $S(f)=\Gamma(s)(f)$. We will proceed as follows : the components of the vector
$\Gamma(s)(f)$ are in $B_0=\Gamma(U_k,\O_{C_k})$, thus their only
poles can occur at the points $\{P_1,\dots,P_r\}$ of $(C \cap \{X_0=0\})(\bar{k})$. In order to estimate these pole orders, we fix a local parameter $t$ at each point $P$ of $(C \cap \{X_0=0\})(\bar{k})$. From above, we can write (in $K_C$ seen as the $\O_{C,P}$-module generated by $1,t^{-1},\dots$) :
$$f=a_{-n}t^{-n}+\dots+a_0+\dots+a_{N-1}t^{N-1}+t^Nh,$$
with the $a_i$ in $R_d$ and $h$ in $\O_{C,P}$. Now we just have to estimate the reduced pole orders of the vectors $S(t^i)$, $S(t^N h)$ ; at the end of the section, we do so for the vectors $S(t^i)$, $i\in \Z$, but we first get rid of the last term, showing that for sufficiently large $N$, the vector $S(t^Nh)$ is in $W_l(\O_{C_k,P})$ for each $h$ in $\O_{C,P}$. 

\medskip

It is sufficient to show that for any $h$ in $\O_{C,P}$, the components of $S(h)$ have bounded pole orders at $P$. Indeed, if $S(t)=(t_0,\dots,t_{l-1})$, the first component $t_0$ is a uniformizing parameter for $C_k$ at $P$, and $v_P(t_0)=1$. For this reason, if we take $N$ sufficiently large, the vector $S(t^Nh)=S(t)^NS(h)$ is in $W_l(\O_{C,P})$.

\bigskip

\medskip

{\bf Lemma 4.2.}
{\it Let $P$ be a point of $C_k$, $n\in \N$. As usual, $v_P$ stands for the valuation at $P$. The set :
$$\RR_{n,P}=\{(f_0,\dots,f_{l-1})\in W_l(K(C_k)),~v_P(f_i)\geq -ip^in\}\subset W_l(K(C_k))$$
is a ring. Moreover, the elements of $\RR_{n,P}$ with $v_P(f_0)=0$ are units of this ring.}

\begin{proof}
The first assertion is straightforward after lemma 0.1. If $f=(f_0,\dots,f_{l-1})$ is in $\RR_{n,P}$ with $v_P(f_0)=0$, let $g=(g_0,\dots,g_{l-1})=f^{-1}$ ; we have :
 $$g=(f_0^{-1},0,\dots,0)+\sum_{i=1}^{l-1} (-1)^iV^i\left(F^i(f_0^{-i-1},0,\dots,0)F^{i-1}(f_1,\dots,f_{l-1},0)^i\right).$$
 Since $f_0$ has valuation $0$, we just have to evaluate the valuations of the other terms. Set $(h_0^{(i)},\dots,h_{l-1}^{(i)}):=(f_1,\dots,f_{l-1},0)^i$. We have $v_P(h_j^{(1)})\geq-(j+1)p^{j+1}n$, and an easy induction on $i$, with the help of lemma 0.1 gives $v_P(h_j^{(i)})\geq-(j+i)p^{j+1}n$. If $(k_0^{(i)},\dots,k_{l-1}^{(i)})=F^{i-1}(f_1,\dots,f_{l-1},0)^i$, we get $v_P(k_j^{(i)})\geq -(j+i)p^{j+i}n$. Since this function appears in the $(i+j)$-th component of $g$, we get $v_P(g_i)\geq -ip^in$, and finally $g \in \RR_{n,P}$. 
\end{proof}

\medskip

We are ready to show the following proposition :

\medskip

{\bf Proposition 4.3.}
{\it Let $D$ and $D_0$ be the divisors defined in theorem 2.1, $P$ a closed point of $C\backslash U$. If $v$ and $v_0$ stand for the multiplicities of $D$ and $D_0$ at $P$, then $S(\O_{C,P})\subset \RR_{v+2v_0,P}$.}

\begin{proof}
Let $P\in \P^n(\bar{k})$, $P\in C(\bar{k})\backslash U(\bar{k})$, and assume $P$ has invertible $n$-th coordinate, that is $P$ is in $U_n(\bar{k})=C\cap \{X_n\neq 0\}(\bar{k})$. Then the regular functions on $U_n$ are the images of the elements of :
$$R[\frac{X_0}{X_n},\dots,\frac{X_{n-1}}{X_n}]=R[\frac{1}{x_n},\dots,\frac{x_{n-1}}{x_n}].$$
Moreover, the elements of $\O_{C,P}$ are the images of the elements of $\O_{\P^n,P}$, that is the $fg^{-1}$, $f,g \in R[\frac{1}{x_n},\dots,\frac{x_{n-1}}{x_n}]$, $\bar{g}(P)\neq 0$ in $\bar{k}$. Note that this last condition is equivalent to asking that the 0-th component of $S(g)$ has valuation zero at $P$. 

Thus if we show that $S(\frac{1}{x_n}),\dots,S(\frac{x_{n-1}}{x_n})$ are in $\RR_{v+2v_0,P}$, the first part of lemma 4.2 ensures us that $S(B_n)\subset \RR_{v+2v_0,P}$, and the second part, joint with the preceding discussion, that $S(\O_{C,P})\subset \RR_{v+2v_0,P}$. 

\medskip

From the description of $S$, we have :
\begin{eqnarray*}
S(\frac{x_i}{x_n})
&
=
&
\delta(\frac{x_i}{x_n})+\sum_{k=1}^{l-1}(-1)^k\delta(x_ix_n^{-k-1})(0,x_{n,1},\dots,x_{n,l-1})^k\\
&
+
&
\sum_{k=1}^{l-1}(-1)^{k-1}\delta(x_n^{-k})(0,x_{n,1},\dots,x_{n,l-1})^{k-1}(0,x_{i,1},\dots,x_{i,l-1}).
\end{eqnarray*}	 
Now the terms $\delta(\frac{x_i}{x_n})$, $\delta(x_ix_n^{-k-1})$ and $\delta(x_n^{-k})$ are in $W_l(\O_{C_k,P})\subset \RR_{v+2v_0,P}$ ; from corollary 2.3, the terms $(0,x_{i,1},\dots,x_{i,l-1})$ are in $\RR_{v+2v_0,P}$, and we are done.
\end{proof}

\bigskip

Now we estimate the reduced pole orders of the vectors $S(t^m)$ ; we first treat the case $P \in C(k)$ ; let $t=x_i/x_n-\alpha,~\alpha\in\T$ be a local parameter for an $R$-point above $P$ as in proposition 4.1. From the description of $S$, we have :
\begin{eqnarray*}
S(t)
&
=
&
(t_0,\dots,t_{l-1})\\
&
=
&
\delta(\frac{x_i}{x_n}-\alpha)+\sum_{h=1}^{l-1}(-1)^h\delta(x_ix_n^{-h-1})(0,x_{n,1},\dots,x_{n,l-1})^h\\
&
+
&
\sum_{h=1}^{l-1}(-1)^{h-1}\delta(x_n^{-h})(0,x_{n,1},\dots,x_{n,l-1})^{h-1}(0,x_{i,1},\dots,x_{i,l-1}).
\end{eqnarray*}

From theorem 2.1, we have : 
$$(0,x_{i1},\dots,x_{i,l-1})=\sum_{k=1}^{l-1}p^{k-1}V\left(X_i^{(k)}-(Fx_i)X_0^{(k)}\right)=\sum_{k=1}^{l-1}V^k\left(F^{k-1}X_i^{(k)}-(F^kx_i)F^{k-1}X_0^{(k)}\right)$$

When we develop the $h$-th and $(h-1)$-th power in the expression for $S(t)$, we get that the general term of the first and second sum are respectively :

$$(-1)^h V^{j_1+\dots+j_h}\left(F^{j_1+\dots+j_h}(x_ix_n^{-h-1}x_n^a)F^{j_1+\dots+j_h-1}(\prod X_n^{(j_i)}\prod_{\#\{j_i\}=a} X_0^{(j_i)})\right)~;$$

$$(-1)^{h-1}V^{j_1+\dots+j_h}\left(F^{j_1+\dots+j_h}(x_i^{\epsilon}x_n^{-h}x_n^{a-\epsilon})F^{j_1+\dots+j_h-1}(\prod X_n^{(j_i)} X_i^{(j_i)^\epsilon}\prod_{\#\{j_i\}=a} X_0^{(j_i)})\right),$$

where $0\leq a \leq h $ is the number of times we choose the term $(Fx_n)X_0^{(j)}$ or $(Fx_i)X_0^{(j)}$ in the development, and $\epsilon$ is $0$ or $1$ depending on whether we choose $F^{k-1}X_i^{(j)}$ or $(F^kx_i)F^{k-1}X_0^{(j)}$ in the last parenthesis. Note that the terms with $a=h$ and $\epsilon=1$ are in one-to-one correspondance in the two sums, with coefficient $(-1)^h$ in the first and $(-1)^{h-1}$ in the second. Thus their sum vanishes and we can impose $a\leq h-1$.

Now the terms of the two sums can be treated similarly regarding to their valuation ; we choose to write the general term of $S(t)$ as :

$$ V^J(F^J(x_ix_n^{-h-1}x_n^a,0,\dots,0)F^{J-1}(\prod_{i=1}^h X_*^{(j_i)})),~J=j_1+\dots+j_h \geq h,~a<h.$$

Remark that $\delta(t)=(x_ix_n^{-1},0,\dots,0)-(a,0,\dots,0):=(t_0,f_1,\dots,f_{l-1})$ is in $W_l(\O_{C_k,P})$, and that $t_0=x_ix_n^{-1}-a$ is a local parameter for $C_k$ at $P$, thus $v_P(t_0)=1$. We can now write the image of $t^m$ :
$$S(t^m)=(t_0^m,0,\dots,0)+\sum_{k=1}^{l-1}(-1)^k c_{m,k} (t_0^{m-k},0,\dots,0)(0,t_1,\dots,t_{l-1})^k,$$
where $c_{m,k}=\binom{-m+k-1}{k}$ if $m<0$, and $c_{m,k}=\binom{m}{k}$ if $m>0$. From above, we can write the general term of $(0,t_1,\dots,t_{l-1})^k$ as :
$$ (0,f_1,\dots,f_{l-1})^{k_0}\left(\prod_{u=1}^t \left(V^{J_u}(F^{J_u}(x_ix_n^{a_u-h_u-1})F^{J_u-1}(\prod_{i_u=1}^{h_u}X_*^{(j_{i_u})}))\right)^{k_u}\right),$$
where $k_0+\dots+k_t=k$. Setting $K=k_0+\sum_{u=1}^t k_uJ_u$, we get :
$$V^{K}\left(F^{K}(t_0^{m-k}\prod_{u=1}^t(x_ix_n^{a_u-h_u-1})^{k_u})F^{K-1}\left(\prod_{u=1}^t(\prod_iX_*^{j_{i_u}})^{k_u}(f_1,\dots,f_{l-1})^{k_0}\right)\right).$$
which is congruent modulo $\wp W_l(K(C_k))$ to :

$$V^{K}\left(F(t_0^{m-k}\prod_{u=1}^t(x_ix_n^{a_u-h_u-1})^{k_u})\prod_{u=1}^t(\prod_iX_*^{j_{i_u}})^{k_u}(f_1,\dots,f_{l-1})^{k_0}\right).$$
We want to evaluate the reduced pole order of this Witt vector. First set $(g_0,\dots,g_{l-1})=\prod_{u=1}^t(\prod_iX_*^{j_{i_u}})^{k_u}$ ; from lemma 0.1, we get that $v_P(g_i)\geq -p^{i+1}(K-k_0)(v+v_0)$ since $\sum k_u\sum j_{i_u}=K-k_0$. 

\smallskip

Next consider the term $H= F(t_0^{m-k}\prod_{u=1}^t(x_ix_n^{a_u-h_u-1})^{k_u})$ ; since $v(x_n)=-v_0$ and $v(x_i)\geq -v_0$, we have $v_P(H)\geq p(m-k-\sum k_u(a_u-h_u)v_0)\geq p(m-k+(k-k_0)v_0)$ since $a_u\leq h_u-1$. 

\smallskip

Setting $(h_0,\dots,h_{l-1})=H(g_0,\dots,g_{l-1})(f_1,\dots,f_{l-1},0)^{k_0}$, and remarking that the last vector has no incidence on our estimation, we get : $v_P(h_i)\geq p^{i+1}(m-k+(k-k_0)v_0-(K-k_0)(v+v_0))$, and :
$$\rp_P(V^K(h_0,\dots,h_{l-1}))\leq -p^{l-K}(m-k+(k-k_0)v_0-(K-k_0)(v+v_0)).$$
Finally, since the $0$-th component of $S(t^m)$ is $t_0^m$, we get $\rp_P(S(t^m))\leq \max(-p^{l-1}m,M)$, where :
\begin{eqnarray*}
M
&
=
&
\max_{1\leq k\leq K\leq l-1~;~0\leq k_0\leq k}\left(-p^{l-K}(m-k+(k-k_0)v_0-(K-k_0)(v+v_0))\right)\\
&
=
&
\max_{1\leq k\leq K\leq l-1~;~0\leq k_0\leq k}\left(-p^{l-K}(m-k+kv_0-K(v+v_0)+k_0v)\right)\\
\end{eqnarray*}

Since $D$ is effective, we have $v\geq 0$, and we have to take $k_0=0$. Now if we fix $K$ and let $k$ vary between $1$ and $K$, the maximum is obtained for $k=1$, and we get $M=\max_{1\leq K\leq l-1}M_K$, $M_K=-p^{l-K} (m-1+v_0-K(v+v_0))$. 
Assume $m<0$ ; computing the $M_K-M_{K+1}$, we get the maximum for $K=1$ or $K=2$, depending on whether $(1-p)(m-1)+v(p-2)-v_0$ is positive or negative, that is :

\medskip

{\bf Lemma 4.4.}
\label{rpS}
{\it Notations being as above, let $m<0$ ; we have :
$$ \rp_P(S(t^m)) \leq \max(p^{l-1}(-m+1+v),p^{l-2}(-m+1+v_0+2v)),$$
and the first term is the maximum if $(1-p)(m-1)+v(p-2)-v_0$ is positive.}

\medskip

In the case $m>0$, we don't need to get results as precise as before : it is sufficient to remark that $\rp_P(S(t^m))\leq \rp_P(S(t))$ here ; by the same calculations as above, we obtain $\rp_P(S(t))\leq \max(p^{l-1}v,p^{l-2}(v_0+2v))$. 

\bigskip

We drop the assumption $P\in C(k)$ : assume now $P\in C(k_d)$ (i.e. $P$ is a point of degree $d$ of the curve $C_k$, via the identification between $C(\bar{k})$ and $C_k(\bar{k})$). Let $t=\pi_d(x_i/x_n)$ be a local parameter for an $R_d$-point above $P$ as in proposition 4.1, and $f \in t^{-n}\O_{C,P}$.

Consider the curve $C_d:=C\times_{\Sp R} \Sp R_d$, its reduction $C_{k_d}:=C_d\otimes k_d$, and let $P_d$ be one of the $d$ points in $C_d(k_d)$ dividing $P$. In $K_{C_d}=K_C\otimes R_d$, we write :
$$\pi_d(\frac{x_i}{x_n})=\prod_{\sigma \in \Gal(R_{d'}/R)} \left( \frac{x_i}{x_n}-\alpha^\sigma\right),$$
for some $d'$ dividing $d$. Exactly one of the $t_d^\sigma=\frac{x_i}{x_n}-\alpha^\sigma$ (say $t_d=\frac{x_i}{x_n}-\alpha$) is a local parameter for an $R_d$-point above $P_d$, and the others are elements of $\O_{C_d,P_d}^*$. Therefore we have $f\in t_d^{-n}\O_{C_d,P_d}$ if we consider $f$ as an element of
$K_{C_d}=K_C\otimes R_d$.

Now it is easy to verify that the morphism $S$ extends to $S_d : K_{C_d}\rightarrow W_l(K(C_{k_d}))$ (it is sufficient to extend $S$ on the constants, i.e. as the Witt isomorphism from $R_d$ to $W_l(k_d)$). If we apply lemma 4.4 to $C_d$ and $P_d$ ($t_d$ is exactly the local parameter given by proposition 4.1 for $C_d$ at $P_d$ since $\alpha\in \T_d$), we get $\rp_{P_d}(S_d(t^m))\leq \max(p^{l-1}(-m+1+v),p^{l-2}(-m+1+v_0+2v)) $. Now the following lemma generalizes lemma 4.4 to points of arbitrary degree :

\medskip

{\bf Lemma 4.5.} 
{\it Let $P \in \P(K)$ a place of degree $d$ and $P_1$ be a place dividing $P$ in $\P(K\otimes k_d)$. If $f\in W_l(K)$ is such that $\rp_{P_1}(f)=a>0$ (in $W_l(K\otimes k_d)$), then $\rp_P(f)=a$ (in $W_l(K)$).}

\medskip

Finally, if we write $f=a_{-n}t^{-n}+\dots+a_{N-1}t^{N-1}+t^Ng$, $ g \in \O_{C,P}$, the preceding results give the :

\medskip

{\bf Proposition 4.6.}
{\it Let $P\in C(\bar{k})\backslash U(\bar{k})$, and $t$ a local parameter at $P$ as in proposition 4.1 ; if $f\in t^{-n}\O_{C,P}$, $n>0$, we have :
$$\rp_P(S(f))\leq \max(p^{l-1}(n+1+v),p^{l-2}(n+1+v_0+2v)).$$
 If $f\in \O_{C,P}$, we get $\rp_P(S(f))\leq \max(p^{l-1}v,p^{l-2}(v_0+2v))$.}

\medskip

{\bf Remark.}
In the case of projective line, if we define the morphism $S:R(x)\rightarrow W_l(k(x))$ as $w_l$ on $R$ and $S(x)=(x,0,\dots,0)$, the example following lemma 1.1, joint with theorem 1.2 show that we are estimating the same character sums as in \cite{li} section 4 ; the calculations we did in this section, joint with the estimation of the degree of the idele class character in section 3, give theorem 4.1 of \cite{li} and its corollaries.

\section{Bounds}

We know come to the estimation of the exponential sums ; let $C\subset \P_R^n$ be a smooth projective curve over $\Sp R$ with $C_k=C\otimes k$ of genus $g$, $H_0$ the hyperplane $\{X_0=0\}\subset \P_R^n$, and $U$ the affine open $C\cap \{X_0\neq 0\}$. Denote by $P_1,\dots,P_r$ the points of $(C\cap H_0)(\bar{k})$ ; let $t_i=\pi_i(x_jx_n^{-1})$ be a local parameter at $\Pi_i$ an $R_d$-point above $P_i$ ($1\leq i\leq r$) as in propostion 4.1, and $\O_C((\Pi_i))$ the associated invertible sheaf. Denote by $D_0$ the intersection divisor of $H_0$ with $C_k$ in $\P_k^n$, and by $v_{0i}$ its multiplicity at $P_i$. Let $D$ be a divisor of degree $2g-2+\lceil \frac{2g-1}{p} \rceil$, whose support is contained in that of $D_0$, and $v_i$ its multiplicity at $P_i$.

\begin{theorem}
\label{fin}
There exists a subset $\T_\sigma \subset U(R)$ which is the image of a section of the map $r_U:U(R)\rightarrow U_k(k)$, and such that for $f=f_0+pf_1+\dots+p^{l-1}f_{l-1}\in \Gamma(U,\O_C)$ a nondegenerate function (that is $S(f)$ is a nondegenerate vector of functions) with
$ f_j\in \Gamma(C,\O_C(\sum_{i=1}^r n_{ij}(\Pi_i))),$
we have the bound :
$$\sum_{\Pi\in \T_\sigma} \exp\left( \frac{2i\pi}{p^l} \Tr(f(\Pi))\right)
\leq \left(\sum_{i=1}^r (A_i+1)\deg P_i+2g-2\right) p^\frac{m}{2},$$
where :
$$A_i=\max\left(\max_{0\leq j\leq l-2}\left(p^{l-1-j}(n_{ij}+1+v_i),p^{l-2-j}(n_{ij}+1+v_{0i}+2v_i)
\right), n_{i,l-1}\right).$$
\end{theorem}

\begin{proof}
First choose for $\T_\sigma$ the Teichm\"uller subset corresponding to a section $\sigma$ satisfying the requirements of theorem 2.1. By the definition of the Cartier divisors $(\Pi_i)$, we must have :
$$\forall i\in \{1\dots,r\}, j\in\{0,\dots,l-1\}, f_j \in t_i^{-n_{i,j}}\O_{C,P_i}.$$

We will get the result from theorems 1.2 and 3.1, once we have estimated reduced pole orders of the Witt vector $\Gamma(s)(f)$ ; first note that :
$$\Gamma(s)(f)
\equiv
\Gamma(s)(f_0)+V\Gamma(s)(f_1)+\dots+V^{l-1}\Gamma(s)(f_{l-1})\quad [\wp W_l(K)]$$
Now applying proposition 4.6, we obtain :
$$\rp_{P_i}(V^j\Gamma(s)(f_j))\leq \max\left(p^{l-j-1}(n_{ij}+1+v_i),p^{l-j-2}(n_{ij}+1+v_{0i}+2v_i)\right),$$
if $j\leq l-2$ ; for $j=l-1$, since $V^{l-1}\Gamma(s)(f_{l-1})=(0,\dots,0,\overline{f_{l-1}})$, the reduced pole order of this vector at $P_i$ is at most $n_{i,l-1}$. 
\end{proof}

We write this result in the particular case $l=2$ : it includes
the Galois ring $\Z/4\Z$ which concentrates most of the efforts
of coding theorists. Moreover the result is simpler : looking
back at the estimations of the reduced pole orders, the maximum
is $p(-m+1+v_i)$ since we can only have $K=1$. The equality
$$\sum_{i=1}^r v_i \deg P_i =2g-2+\lceil\frac{2g-1}{p}\rceil$$
(this is the degree of the divisor $D$) allows us to give a bound depending only on the pole orders of the function and on the genus of the curve $C_k$ :

\medskip

{\bf Corollary 5.2.}
\label{fin1}
{\it Let $R=GR(p^2,m)$ be a Galois ring, and $C$ a smooth curve over $R$, $C\subset \P_R^n$. We can find a subset $\T_\sigma \subset U(R)$ which is the image of a section of the map $r_U:U(R)\rightarrow U_k(k)$, and such that for $f\in \Gamma(C,\O_C(\sum n_i(\Pi_i)))$ a nondegenerate function, we have the bound :
$$\sum_{\Pi\in \T_\sigma} \exp\left( \frac{2i\pi}{p^2} \Tr(f(\Pi))\right)
\leq B_f p^{\frac{m}{2}},$$
where }
$$B_f=\sum_{i=1}^r (p(n_i+1)+1)\deg P_i +(2g-2)(p+1)+p\lceil
\frac{2g-1}{p}\rceil.$$

\medskip
It can be shown that, if $C_k\subset \P_k^n$ is a smooth curve, either a complete intersection or a canonical curve (that is, the image of a curve $C$ of genus $g$ by the canonical embedding $C\subset \P_k^{g-1}$), then it lifts to a smooth curve $C\subset \P_R^n$. In this case, the degree of the divisor $D$ is greater than that
of $D_0$, we can choose $v_i \geq v_{0i}$ for each
$i\in\{1,\dots,r\}$ ; if moreover we assume $p\geq 3$, the maximum in proposition 4.6 is the first term, and as above the bound does not depend on the projective embedding :

\medskip

{\bf Corollary 5.3.}
\label{fin2}
{\it Let $C$ be the lift over $\Sp R$ of either a canonical or a complete intersection
curve. Assume $p\neq 2$.  We can find a subset $\T_\sigma \subset U(R)$ which is the image of a section of the map $r_U:U(R)\rightarrow U_k(k)$, and such that for $f\in \Gamma(C,\O_C(\sum n_i(\Pi_i)))$ a nondegenerate function,
we have the bound :
$$\sum_{\Pi\in \T_\sigma} \exp\left( \frac{2i\pi}{p^l} \Tr(f(\Pi))\right)
\leq B_f p^{\frac{m}{2}},$$
where }
$$B_f=  \sum_{i=1}^r (p^{l-1}(n_i+1)+1)\deg P_i
+(2g-2)(p^{l-1}+1)+p^{l-1}
\lceil\frac{2g-1}{p}\rceil $$

\end{document}